\documentclass[preprint,12pt,authoryear]{elsarticle}

\usepackage{amssymb,amsmath,amsthm,amsfonts}
\usepackage[linesnumbered,ruled,noend]{algorithm2e}
\usepackage{textcomp}
\usepackage{xargs}
\usepackage[svgnames]{xcolor}
\usepackage{hyperref}
\hypersetup{
    colorlinks,
    citecolor=CornflowerBlue,
    filecolor=CornflowerBlue,
    linkcolor=RoyalBlue,
    urlcolor=CornflowerBlue
}
\definecolor{purpleweb}{rgb}{0.5,0.0,0.5}
\usepackage{todonotes}
\newcommandx{\iman}[2][1=]{\todo[linecolor=CornflowerBlue,backgroundcolor=CornflowerBlue!25,bordercolor=CornflowerBlue,size=small,author=Iman,#1]{#2}}
\newcommandx{\elad}[2][1=]{\todo[linecolor=purpleweb,backgroundcolor=purpleweb!25,bordercolor=purpleweb,size=small,author=Elad,#1]{#2}} 
\usepackage{cleveref}
\usepackage{mathtools}
\usepackage{setspace}

\DeclareMathOperator*{\lexmax}{lexmax}

\DeclareMathOperator*{\argmax}{argmax}
\DeclareMathOperator*{\argmin}{argmin}
\newtheorem{theorem}{Theorem}
\newtheorem{lemma}{Lemma}

\newtheorem{remark}{Remark}
\newtheorem{assumption}{Assumption}
\newtheorem{definition}{Definition}

\newtheorem{example}{Example}
\newtheorem{problem}{Problem}
\SetKwBlock{Repeat}{Forever}{}

\journal{European Journal of Operational Research}

\begin{document}

\begin{frontmatter}



\title{Sensitivity Analysis for Bottleneck Assignment Problems}


\author[1]{Elad~Michael\corref{cor1}}\ead{eladm@student.unimelb.edu.au}
\author[2]{Tony~A.~Wood}\ead{tony.wood@epfl.ch}
\author[1]{Chris~Manzie}\ead{manziec@unimelb.edu.au}
\author[3]{Iman~Shames}\ead{iman.shames@anu.edu.au}

\cortext[cor1]{Corresponding author}
\address[1]{Department of Electrical and Electronic Engineering, University of Melbourne, Parkville, Victoria, 3010, Australia}
\address[2]{SYCAMORE Lab, École Polytechnique Fédérale de Lausanne (EPFL), Lausanne, Switzerland}
\address[3]{CIICADA Lab, School of Engineering, Australian National University, Acton, ACT, 0200, Australia}

\begin{abstract}
In assignment problems, decision makers are often interested in not only the optimal assignment, but also the sensitivity of the optimal assignment to perturbations in the assignment weights. Typically, only perturbations to individual assignment weights are considered. We present a novel extension of the traditional sensitivity analysis by allowing for simultaneous variations in all assignment weights. Focusing on the bottleneck assignment problem, we provide two different methods of quantifying the sensitivity of the optimal assignment, and present algorithms for each. Numerical examples as well as a discussion of the complexity for all algorithms are provided.
\end{abstract}


\end{frontmatter}


\section{Introduction}\label{sec:Intro}
In the classic assignment problem, we seek a one to one matching of a set of agents to a set of tasks which optimises some assignment cost. For example, \cite{nam2015your} consider a vehicle routing problem, assigning destinations to vehicles in order to minimise the time taken for all vehicles to reach their destinations within a road network. Assignment problems have a wide variety of applications, such as resource allocation (\cite{harchol1999choosing}), scheduling (\cite{carraresi1984multi}, \cite{adams1988shifting}), and the aforementioned autonomous vehicle routing (\cite{8264324}, \cite{nam2015your}, \cite{arslan2007autonomous}). \cite{pentico2007assignment} review a variety of assignment objectives, such as the sum of all the assigned weights (linear assignment) or the difference between the maximum and minimum weight assignment (balanced assignment problem). Here we focus on the bottleneck assignment problem (BAP), which minimises the maximum assigned weight within the assignment. The BAP is particularly applicable in minimum time scenarios, ensuring that a team of agents working in parallel complete their tasks in minimum time. For the bottleneck assignment problem there are many algorithms which provide solution in polynomial time, reviewed in~\cite{pentico2007assignment} and~\cite{pferschy1997solution}.

In many applications the weight associated with the assignment of an agent to a task carries with it some uncertainty. Measurement errors, communication delays, or a dynamic environment may contribute to an inaccurate estimate of an assignment weight. If the uncertainty is known \textit{a priori} a \emph{minimax} approach may be used, identifying an assignment which minimises the assignment cost for the worst case realisation of the uncertainty. However, this is typically a conservative solution. In applications where such \textit{a priori} knowledge is unavailable/impractical, or a conservative solution undesirable, we instead employ sensitivity analysis.

Sensitivity analysis characterises how perturbations to the inputs affect the output. In the case of uncertainty in the assignment weights, a sensitivity analysis informs the decision maker in which ways the optimal assignment may change, given perturbations to the measured weights. We adopt Type II sensitivity analysis from \cite{koltai2000difference}, or equivalently the sensitivity of the optimal basis as in \cite{jansen1997sensitivity}. Type II sensitivity analysis characterises the set of weight perturbations to which the optimal assignment is invariant. In the sensitivity analysis of linear programs, as in \cite{jansen1997sensitivity}, the perturbations are restricted to an individual element of the canonical weight vector, while the remaining elements are held fixed. However, due to the additional structure of a BAP, we are able to extend the sensitivity analysis to simultaneous perturbations in all of the weights. We consider this a more applicable formulation of assignment sensitivity, as perturbations to the measured weights occur simultaneously if tasks are carried out in parallel.

\subsection{Literature Review}\label{sec:litReview}

Due to the assignment problem's popularity in network theory, there exist various methods incorporating uncertainty in the assignment weights. Roughly characterised, we have those which compute robust assignments with \textit{a priori} knowledge of uncertainty, and those which characterise some notion of sensitivity of the assignment to perturbations. \cite{volgenant2010improved} presents several robust optimisation approaches to the bottleneck assignment problem, with interval assumptions on the possible weight values. \cite{fu2015robust} and \cite{kasperski2013bottleneck} begin with a set of possible scenarios, and compute robust minimax assignment solutions by blending those scenarios. Without \textit{a priori} knowledge of the uncertainty, \cite{lin2003sensitivity} and \cite{volgenant2006addendum} present the sensitivity of the linear assignment to perturbations in a single weight, based on the dual variables arising from a primal-dual optimisation. \cite{ramaswamy2005sensitivity} provide similar results in the shortest path and maximum cut network problems, which can be shown to be equivalent formulations to the linear assignment problem. However, they do not consider simultaneous weight perturbations in multiple edges.
Both \cite{nam2015your} and \cite{lin2007sensitivity} present sensitivity analysis algorithms for the linear assignment, allowing for perturbations in an assigned edge as well as in edges adjacent to the assigned edge (i.e., a row/column of the weight matrix). However, these analyses do not take into account simultaneous perturbations in multiple assigned edges. \cite{wood2020collision} present an assignment and collision avoidance algorithm, which utilises a bottleneck assignment sensitivity estimation. However, they consider only positive weight perturbations, and only in the edges which are assigned.

Most similar to the analysis presented in this paper are the results from \cite{sotskov1995some}, which discusses the largest \emph{uniform} bound on all of the input weight perturbations to maintain optimality for both the linear and bottleneck assignment problems. The results from \cite{sotskov1995some} will coincide with the smallest magnitude perturbation bound in our results, discussed in Section~\ref{sec:probForm} and Section~\ref{sec:complexity}.

Previously, \cite{michael2019uncertainty} presented a preliminary sensitivity analysis which allowed for simultaneous perturbations in all assignment weights. In this paper we build on these results, providing a complete theoretical framework for understanding the sensitivity of the bottleneck assignment problem, as well as improved algorithmic results. For a given bottleneck assignment, the methods we provide construct intervals within which the assignment weights may vary while preserving the optimality of the assignment. We further show these are the \emph{lexicographically largest} set of allowable perturbations, as defined in Section~\ref{sec:probForm}. We provide algorithms to compute these intervals, as well as numerical examples and complexity analyses.

\subsection{Organisation of Paper}
Section~\ref{sec:probForm} formalises the assignment problem, as well as the definitions of sensitivity used in this work. Section~\ref{sec:Esensitivity} establishes the main theoretical tools used, and applies them to solve the problem of bottleneck \emph{edge} sensitivity. Using these, the bottleneck assignment sensitivity problem is addressed in Section~\ref{sec:sensBttleAss}.
In Section~\ref{sec:complexity}, we discuss the complexities of the provided algorithms. Conclusions are presented in Section~\ref{sec:conclusion}.
\section{Preliminaries and Problem Formulation}\label{sec:probForm}

Let $\mathcal{G}=(\mathcal{V}, \mathcal{E})$ be a bipartite graph with $\mathcal{V} = \mathcal{V}_{1} \cup \mathcal{V}_{2}$ the set of vertices such that $\mathcal{V}_{1} \cap \mathcal{V}_{2} = \varnothing$, and the edge set $\mathcal{E} \subseteq \mathcal{V}_{1} \times \mathcal{V}_{2}$. While formulating the assignment problem, it is convenient to use the pair of vertices that an edge connects to represent the edge, i.e. edge $(i,j)\in\mathcal{E}$ from vertex $i$ to $j$. Define a weight matrix $W\in\overline{\mathbb{R}}^{n\times m}$ over the extended reals i.e., $\overline{\mathbb{R}}:= \mathbb{R}\cup\{-\infty,+\infty\}$, where $n:=|\mathcal{V}_{1}|$ and $m:=|\mathcal{V}_{2}|$ such that $n \geq m \geq 1$, and $w_{ij}$ is the weight associated with edge $(i,j)$, or $+\infty$ if $(i,j)\not\in\mathcal{E}$.
We also define a set of binary decision variables $\pi := \{\pi_{ij} \mid \pi_{ij} \in \{0,1\},(i,j) \in \mathcal{E}\}$. The bottleneck assignment problem can be formulated as
\begin{subequations}\label{eq:bttleNeckOpt}
\begin{align}
\underset{\pi_{ij}\in\{0,1\}}{\min}\; \underset{(i,j) \in \mathcal{E}}{\max}\; & \pi _{ij} w_{ij} \label{eq:bttleneckCost}\\
\textrm{subject to } &\sum_{i \in \mathcal{V}_{1}} \pi_{ij} = 1, \; j \in \mathcal{V}_{2} \label{eq:fewerAgentsi}\\
&\sum_{j \in \mathcal{V}_{2}} \pi_{ij} \leq 1, \; i \in \mathcal{V}_{1}. \label{eq:moreTasksj}
\end{align}
\end{subequations}

If $\pi_{ij} = 1$ then we say that vertex $i$ is \emph{assigned} to vertex $j$ or the edge $(i,j)$ is assigned. Constraints \eqref{eq:fewerAgentsi} and \eqref{eq:moreTasksj} ensure that every vertex is assigned to at most one other vertex, and that all vertices from the smaller vertex set $\mathcal{V}_2$ are assigned.
Along with the set of binary decision variables $\pi$, we define an assignment $\Pi = \{(i,j) \mid \pi_{ij} = 1\}$ to be the set of edges which are assigned in $\pi$. We refer to an assignment which is an optimiser of~\eqref{eq:bttleNeckOpt} as a \emph{bottleneck assignment}. In the remainder of the paper, we use $e\in\mathcal{E}$ to represent an edge rather than the pair of vertices $(i,j)$.

\begin{definition}[Feasible Assignments Set]\label{def:assignmentSpace}
For a bipartite graph $\mathcal{G}$ with weights $W\in\overline{\mathbb{R}}^{n\times m}$, define $\mathcal{P}^{n\times m}$ as the set of all possible assignments $\Pi$ satisfying the constraints~\eqref{eq:fewerAgentsi}-\eqref{eq:moreTasksj}. Additionally, define $\mathcal{P}_{e}\subseteq\mathcal{P}^{n\times m}$ to be the subset of assignments such that $e\in\Pi$ for all $\Pi\in\mathcal{P}_{e}$.
\end{definition}
We define a mapping ${\bf BAP}: \overline{\mathbb{R}}^{n\times m} \rightarrow \mathcal{P}^{n\times m}$ which maps the weights to the set of optimisers of~\eqref{eq:bttleNeckOpt}. We briefly describe a simple algorithm which solves the bottleneck assignment problem: Let $\Pi$ be an assignment over the graph. If there does not exist an assignment consisting of edges with strictly lower weight than $\Pi$, then $\Pi$ is optimal. If there does, replace $\Pi$ with the found assignment. Repeat until $\Pi$ is optimal.

\begin{definition}[Allowable Perturbation]\label{def:allowable}
For a given bipartite graph $\mathcal{G}$ with weights $W$, let $\Pi\in{\bf BAP}(W)$, i.e. $\Pi$ is an optimiser of~\eqref{eq:bttleNeckOpt} over the weights $W$. A perturbation $P\in\overline{\mathbb{R}}^{n\times m}$ is \emph{allowable} with respect to $\Pi$ if $\Pi \in {\bf BAP}(W+P)$. If $\Pi \not\in {\bf BAP}(W+P)$, then $P$ is not allowable with respect to $\Pi$.
\end{definition}

For a given graph $\mathcal{G}$ and weights $W$ with optimiser $\Pi\in{\bf BAP}(W)$, identifying if $\Pi$ is an optimiser over the perturbed weights $W+P$ is, in general, equivalent to solving a new assignment problem. In this sensitivity analysis, we provide a sufficient but not necessary condition for the invariance of the optimal assignment, in the form of an interval test.

Let $\Lambda\subseteq\overline{\mathbb{R}}^{n\times m}$ be an $n\times m$ array of intervals over the extended reals. For each edge $e\in\mathcal{E}$, let $[-\underline{\lambda}_{e},\overline{\lambda}_{e}]$ be the interval corresponding to edge $e$. 
\begin{remark}\label{rem:arrayIndex}
As $\Lambda = ([-\underline{\lambda}_{e},\overline{\lambda}_{e}])_{e\in\mathcal{E}}$, we may access the upper and lower bound values $\underline{\lambda}_{e},\overline{\lambda}_{e}$ by their edge indices.
\end{remark}
For any perturbation $P\in\overline{\mathbb{R}}^{n\times m}$, we write $P\in\Lambda$ if
\begin{align*}
{\forall e \in \mathcal{E}:}\; {P[e]} \in [-\underline{\lambda}_{e},\overline{\lambda}_{e}],
\end{align*}
where we have used $P[e]$ to denote the perturbation element corresponding to edge $e$. Finally, we label $\Lambda$ allowable relative to $\Pi$ if for all $P\in\Lambda$, $P$ is an allowable perturbation relative to $\Pi$ as in Definition~\ref{def:allowable}. Throughout this paper we will require that ${\bf 0}^{n \times m}\in\Lambda$, or equivalently $\underline{\lambda}_{e},\overline{\lambda}_{e}$ are non-negative. To order the arrays intervals we use the \emph{lexicographic ordering}.

The \emph{lexicographic} ordering is inspired by alphabetical ordering of words, such as in a dictionary. Let $\Lambda\subseteq\overline{\mathbb{R}}^{n\times m}$ be an array of intervals. Define a mapping $\rho_k(\Lambda)$ which returns the $k$-th minimum magnitude upper or lower bound of any of the intervals in $\Lambda$, which will be used to order these arrays.

\begin{definition}[Lexicographic Ordering of Interval Arrays]
For two arrays of intervals, $\Lambda_a,\Lambda_b\subseteq\overline{\mathbb{R}}^{n\times m}$, $\Lambda_a \succ \Lambda_b$ if there exists $i\in[1,...,nm]$ such that
\begin{align}
\rho_i(\Lambda_a) > \rho_i(\Lambda_b), \textrm{  \emph{and}}\quad {\forall k<i:\;}\rho_k(\Lambda_a) = \rho_k(\Lambda_b).
\end{align}
Let the operator $\lexmax \left( \mathcal{I} \right)$ return a \emph{lexicographic maximal element} of $\mathcal{I}$ denoted by $\Lambda^*\in\mathcal{I}$ where $\mathcal{I}$ is a set of arrays of intervals and $\Lambda^*\succeq \Lambda$, for all $\Lambda\in\mathcal{I}$, with $\succeq$ denoting the lexicographic ordering.
\end{definition}

Lexicographic maximisation thus maximises the minimum element, then the second minimum, and so on.

\begin{problem}\label{prob:BAPSensitivity}
For a given graph $\mathcal{G}=(\mathcal{V},\mathcal{E})$ with weights $W$ and optimiser $\Pi\in{\bf BAP}(W)$, construct the lexicographically largest array of intervals $\Lambda^* := ([-\underline{\lambda}^*_{e},\overline{\lambda}^*_{e}])_{\forall e\in\mathcal{E}}$ such that $\forall e\in\mathcal{E}\;: \underline{\lambda}^*_{e},\overline{\lambda}^*_{e} \geq 0$ and any perturbation $P\in\Lambda^*$ is allowable with respect to $\Pi$.
\end{problem}

For a given graph $\mathcal{G}$ with weights $W$ and $\Pi\in{\bf BAP}(W)$, let
\begin{align}
\mathcal{L}_0 &:= \{\Lambda\subseteq\overline{\mathbb{R}}^{n\times m} \mid {\bf 0}^{n \times m}\in\Lambda\}, \label{eq:L0}\\
\mathcal{L}_A &:= \{\Lambda\subseteq\overline{\mathbb{R}}^{n\times m} \mid {\forall P\in\Lambda:\;} \Pi\in{\bf BAP}(W+P)\}. \label{eq:bapAllowableSet}
\end{align}
Problem~\ref{prob:BAPSensitivity} is then equivalent to finding $\Lambda^*\subset\overline{\mathbb{R}}^{n\times m}$ such that
\begin{align}\label{eq:bapOptOriginal}
\Lambda^* = \lexmax \left( \mathcal{L}_0 \cap \mathcal{L}_{A} \right).
\end{align}
\begin{remark}
Note that the set $\mathcal{L}_0$ enforces the constraint ${\bf 0}^{n \times m}\in\Lambda$, as otherwise the array of intervals could be shifted uniformly by a constant to become arbitrarily ``large''.
\end{remark}
Problem~\ref{prob:BAPSensitivity} is the formal statement of the extension of previous work in sensitivity analysis, as discussed in Section~\ref{sec:litReview}. Note that in general, a lexicographically largest array of intervals will not be unique. However, we make assumptions in the following sections which guarantee that the lexicographic maximum is unique. By constructing the array of intervals $\Lambda^*$, we find the lexicographically largest bounds within which the edge weights may vary while the given assignment $\Pi\in{\bf BAP}(W)$ remains an optimiser. However, the bottleneck assignment $\Pi$ is quite often not unique. The bottleneck assignment problem is concerned only with the maximum weight edge in an assignment, with no mechanism to differentiate between various assignments with the same maximum weight edge. Motivated by this, we may define a related sensitivity problem, focusing only on the maximum weight edge rather than the entire assignment.

\begin{definition}[Bottleneck Edges]\label{def:bottleneckEdge}
For a bipartite graph $\mathcal{G}$ with weights $W$, the maximum weight edges in all bottleneck assignments are the \emph{bottleneck edges} of $\mathcal{G}$.
\end{definition}

We define a new mapping ${\bf E}:\overline{\mathbb{R}}^{n\times m} \rightarrow \mathcal{E}$ which takes the weights and returns the bottleneck edges. We define a perturbation $P\in\overline{\mathbb{R}}^{n\times m}$ as \emph{edge allowable} relative to a bottleneck edge $e^*\in{\bf E}(W)$  if $e^*\in{\bf E}(W+P)$. We formulate the problem of the bottleneck edge sensitivity similarly to the bottleneck assignment sensitivity.

\begin{problem}\label{prob:ESensitivity}
For a given bipartite graph $\mathcal{G}=(\mathcal{V},\mathcal{E})$ with weights $W$ and $e^*\in{\bf E}(W)$, construct the lexicographically largest array of intervals $\Lambda^* := ([-\underline{\lambda}^*_{e},\overline{\lambda}^*_{e}])_{\forall e\in\mathcal{E}}$ such that $\forall e\in\mathcal{E}\;: \underline{\lambda}^*_{e},\overline{\lambda}^*_{e} \geq 0$ and any perturbation $P\in\Lambda^*$ is edge allowable with respect to $e^*$.
\end{problem}
For a bipartite graph $\mathcal{G}$ with weights $W$ and $e^*\in{\bf E}(W)$, let
\begin{align}
\mathcal{L}_{E} &:= \{\Lambda\subseteq\overline{\mathbb{R}}^{n\times m} \mid {\forall P\in\Lambda:\;} e^*\in{\bf E}(W+P)\}. \label{eq:eAllowableSet}
\end{align}
Problem~\ref{prob:ESensitivity} is then equivalent to finding $\Lambda^*\subset\overline{\mathbb{R}}^{n\times m}$ such that
\begin{align}\label{eq:eOptOriginal}
\Lambda^* = \lexmax \left(  \mathcal{L}_0 \cap \mathcal{L}_{E} \right).
\end{align}

This problem will be addressed in Section~\ref{sec:Esensitivity}. We discuss Problem~\ref{prob:ESensitivity} first as the theoretical tools are better motivated and easier to derive in the context of bottleneck edge sensitivity, and will be central in the analysis of bottleneck assignment sensitivity. The problem of assessing the sensitivity of a bottleneck edge is directly comparable to the sensitivity of the bottleneck assignment for the \emph{lexicographic} assignment, see \cite{burkard1991lexicographic}, and the solution of Problem~\ref{prob:ESensitivity} is a conservative estimate of the sensitivity from Problem~\ref{prob:BAPSensitivity}.
\subsection{Exclusive Sets}~\label{sec:exclusive}
We begin by defining the central theoretical tool of this paper, which we term the exclusive set with respect to some $e\in\mathcal{E}$, denoted $\mathcal{S}_{e}$.

\begin{definition}[Exclusive Set]\label{def:exclusiveSet}
For a bipartite graph $\mathcal{G}$ we define $\mathcal{S}_{e^*}$ to be an \emph{exclusive set} with respect to an edge $e^*\in\mathcal{E}$ if it satisfies the following properties:
\begin{enumerate}
    \item $\mathcal{S}_{e^*}\subset\mathcal{E}$;
    \item $e^*\not\in\mathcal{S}_{e^*}$;
    \item $\mathcal{P}^{n\times m} = \mathcal{P}_{e^*} \cup_{e\in\mathcal{S}_{e^*}} \mathcal{P}_{e}$ where $\mathcal{P}_{e}$ is defined in Definition~\ref{def:assignmentSpace}.
\end{enumerate}
For a bipartite graph $\mathcal{G}$, define ${\bf S}_{e}$ to be the set of all exclusive sets for $e\in\mathcal{E}$.
\end{definition}
The third property of an exclusive set $\mathcal{S}_{e^*}$, as defined above, is critical in connecting the exclusive set to the bottleneck assignment problem. In words, we may rephrase the property as: \textit{All assignments over the graph either contain the edge $e^*$, or an edge from the exclusive set}. The connection between the bottleneck assignment problem and the exclusive set is made clear in Lemma~\ref{lem:exclusive}.

\begin{lemma}\label{lem:exclusive}
For a bipartite graph $\mathcal{G}$ with weights $W$, let $e^*\in\mathcal{E}$. There exists an assignment $\Pi\in\mathcal{P}_{e^*}$ and an exclusive set $\mathcal{S}_{e^*}\in{\bf S}_{e^*}$ that satisfy
\begin{align}
w_{e^*} &= \max_{e\in\Pi} w_e \label{eq:maxWeightEdge}\\
w_{e^*} &\leq \min_{e\in\mathcal{S}_{e^*}} w_e \label{eq:lessThanExclusive}
\end{align}
 if and only if $e^*\in{\bf E}(W)$, i.e., the edge $e^*$ is a bottleneck edge, and $\Pi$ is a bottleneck assignment.
\end{lemma}
\begin{proof}
First, assume that edge $e^*\in{\bf E}(W)$ is a bottleneck edge of the graph $\mathcal{G}$ with associated weight matrix $W$. Then there exists at least one assignment $\Pi$ where $e^*$ is the maximum weight edge, so $\Pi$ satisfies~\eqref{eq:maxWeightEdge}. Because $e^*$ is a bottleneck edge of the graph, we have that there is no assignment in which all edges have weight less than $w_{e^*}$. Therefore, all assignments which do not assign $e^*$ contain at least one edge with weight greater than or equal to $w_{e^*}$. Let $\mathcal{S}_{e^*}$ be the set of all edges with weight greater than or equal to $w_{e^*}$, excluding $e^*$. Thus, from Definition~\ref{def:exclusiveSet}, $\mathcal{S}_{e^*}$ is an exclusive set corresponding to $e^*$ and satisfies \eqref{eq:lessThanExclusive}.

To prove the converse, let the assignment $\Pi$ and the edge set $\mathcal{S}_{e^*}$ satisfy~\eqref{eq:maxWeightEdge}--\eqref{eq:lessThanExclusive}, respectively. To obtain a contradiction, we assume there exists an assignment $\Pi'$ that does not assign $e^*$ and only assigns edges with strictly lower weights. We have from the definition of an exclusive set that $\mathcal{P}^{n\times m}\setminus\mathcal{P}_{e^*} \subseteq \cup_{e\in\mathcal{S}_{e^*}} \mathcal{P}_{e}$, i.e., all assignments which do not include $e^*$, assign at least one edge from $\mathcal{S}_{e^*}$. Therefore, there is at least one edge within $\mathcal{S}_{e^*}$ which is assigned in $\Pi'$. However, as $\mathcal{S}_{e^*}$ satisfies~\eqref{eq:lessThanExclusive}, we conclude that $\Pi'$ must contain an edge with weight greater than the weight of $e^*$, which is a contradiction.
\end{proof}

\begin{remark}
The existence of an exclusive set $\mathcal{S}_{e^*}$ satisfying the inequality~\eqref{eq:lessThanExclusive} may be interpreted as proof that there does not exist any augmenting path avoiding $e^*$ with strictly lower costs, see~\cite{jonker1986improving} for more detail on augmenting paths. There may be paths with equivalent costs, as the inequality in~\eqref{eq:lessThanExclusive} is not strict.
\end{remark}

Lemma~\ref{lem:exclusive} establishes the equivalence of $e^*$ being a bottleneck edge and the existence of an exclusive set and assignment $\Pi$ satisfying~\eqref{eq:maxWeightEdge} and \eqref{eq:lessThanExclusive} relative to $e^*$.
This tool allows a reformulation of \eqref{eq:bapOptOriginal} and \eqref{eq:eOptOriginal}, in which we replace the sets $\mathcal{L}_{A}$ and $\mathcal{L}_{E}$ defined in~\eqref{eq:bapAllowableSet} and~\eqref{eq:eAllowableSet}, respectively, with sets redefined by the conditions from Lemma~\ref{lem:exclusive}.

\section{Bottleneck Edge Sensitivity}~\label{sec:Esensitivity}
We now apply Lemma~\ref{lem:exclusive} to reformulate and solve Problem~\ref{prob:ESensitivity}. For a bipartite graph $\mathcal{G}$ with weights $W$, let $e^*\in{\bf E}(W)$.
By Lemma~\ref{lem:exclusive}, there exist an assignment $\Pi$ and an exclusive set $\mathcal{S}_{e^*}$ satisfying~\eqref{eq:maxWeightEdge} and \eqref{eq:lessThanExclusive}, respectively, relative to $e^*$. We define the following sets of interval arrays:
\begin{align}
\mathcal{A}_{\Pi} &:= \{\Lambda \subset \overline{\mathbb{R}}^{n\times m} \mid {\forall e\in\Pi:\;} w_{e^*}-\underline{\lambda}_{e^*} \geq w_{e} + \overline{\lambda}_{e}\}, \label{eq:lSet4Pi} \\
\mathcal{X}_{\mathcal{S}_{e^*}} &:= \{\Lambda \subset \overline{\mathbb{R}}^{n\times m}  \mid {\forall e\in\mathcal{S}_{e^*}:\;} w_{e^*}+\overline{\lambda}_{e^*} \leq w_{e} - \underline{\lambda}_{e}\}. \label{eq:lSet4Exclusive}
\end{align}

The definition of sets $\mathcal{A}_{\Pi},\mathcal{X}_{\mathcal{S}_{e^*}}$ is motivated by Lemma~\ref{lem:exclusive}. Note that for all perturbations in all arrays of intervals of $\mathcal{A}_{\Pi}$, the assignment $\Pi$ has greatest weight edge $e^*$. Similarly, for all perturbations in all arrays of intervals of $\mathcal{X}_{\mathcal{S}_{e^*}}$, we have that $e^*$ has lower weight than any of the edges in $\mathcal{S}_{e^*}$. The main result of this section is presented below.
\begin{theorem}\label{thm:eReformulation}
For a given bipartite graph $\mathcal{G}$ with weights $W\in\overline{\mathbb{R}}^{n\times m}$, let $e^*\in{\bf E}(W)$. There exists an assignment $\Pi\in\mathcal{P}_{e^*}$ and an exclusive set $\mathcal{S}_{e^*}$ satisfying~\eqref{eq:maxWeightEdge} and \eqref{eq:lessThanExclusive}, respectively, such that
\begin{align}
\lexmax(\mathcal{L}_{0}\cap\mathcal{X}_{\mathcal{S}_{e^*}}\cap\mathcal{A}_{\Pi}) = \lexmax(\mathcal{L}_0 \cap \mathcal{L}_{E}) \label{eq:eEquivalence}
\end{align}
where $\mathcal{L}_0$, $\mathcal{L}_{E}$, $\mathcal{A}_{\Pi}$, $\mathcal{X}_{\mathcal{S}_{e^*}}$ are defined in~\eqref{eq:L0}, \eqref{eq:eAllowableSet},\eqref{eq:lSet4Pi}, and \eqref{eq:lSet4Exclusive}, respectively.
\end{theorem}

\begin{proof}
Let $\Lambda^* = \lexmax(\mathcal{L}_0 \cap \mathcal{L}_{E})$. We first prove there exist an assignment $\Pi$ and an exclusive set $\mathcal{S}_{e^*}$ such that $\Lambda^*\in(\mathcal{X}_{\mathcal{S}_{e^*}}\cap\mathcal{A}_{\Pi})$. We then show that $(\mathcal{X}_{\mathcal{S}_{e^*}}\cap\mathcal{A}_{\Pi})\subseteq\mathcal{L}_{E}$, completing the proof.

Recall $\Lambda^* = ([-\underline{\lambda}^{*}_{e},\overline{\lambda}^{*}_{e}])_{e\in\mathcal{E}}$ is an array of intervals, corresponding to each edge in $\mathcal{E}$. For $P[e]$ the perturbation corresponding to edge $e$, we select two perturbations $P_1,P_2\in\Lambda^*$ with elements
\begin{align*}
P_1[e] :=
\begin{cases}
-\underline{\lambda}^{*}_{e} &\textrm{if} \; e \neq e^*\\
\overline{\lambda}^{*}_{e} &\textrm{if} \; e = e^*
\end{cases},\qquad
P_2[e] :=
\begin{cases}
\overline{\lambda}^{*}_{e} &\textrm{if} \; e \neq e^*\\
-\underline{\lambda}^{*}_{e} &\textrm{if} \; e = e^*
\end{cases}.
\end{align*}
Every edge perturbation in $P_1$ is the lower bound of its corresponding interval in $\Lambda^*$, except $e^*$ which is at the upper bound, and the opposite is true for $P_2$. By~\eqref{eq:eAllowableSet}, we have that $e^*\in{\bf E}(W+P)$, for all $P\in\Lambda^*$, including perturbations $P_1$ and $P_2$.
By Lemma~\ref{lem:exclusive}, we have that there exists an assignment and an exclusive set satisfying~\eqref{eq:maxWeightEdge} and \eqref{eq:lessThanExclusive} for each of the perturbations $P_1$ and $P_2$.
We let $\mathcal{S}_{e^*}$ be defined as the exclusive set satisfying~\eqref{eq:lessThanExclusive} for $W_1 := W+P_1$, and $\Pi$ be defined as the assignment satisfying~\eqref{eq:maxWeightEdge} for $W_2 := W+P_2$. By the construction of $P_1$, we have that
\begin{align*}
{\forall e\in\mathcal{S}_{e^*},\;\forall P\in\Lambda^*:\;} w_{e^*}+P[e^*] \leq w_{e^*} + \overline{\lambda}^{*}_{e^*} \leq w_{e} - \underline{\lambda}^{*}_{e} \leq w_{e} - P[e].
\end{align*}
Similarly for $P_2$,
\begin{align*}
{\forall e\in\Pi\setminus\{e^*\},\;\forall P\in\Lambda^*:\;} w_{e}+P[e] \leq w_{e} + \overline{\lambda}^{*}_{e} \leq w_{e^*} - \underline{\lambda}^{*}_{e^*} \leq w_{e^*} - P[e^*].
\end{align*}
Therefore, \eqref{eq:maxWeightEdge} and \eqref{eq:lessThanExclusive} are satisfied by $\Pi$ and $\mathcal{S}_{e^*}$ for all $P\in\Lambda^*$, and thus $\Lambda^* \in \mathcal{A}_{\Pi}$, $\Lambda^* \in \mathcal{X}_{\mathcal{S}_{e^*}}$, and
\begin{align*}
\Lambda^* = \lexmax(\mathcal{L}_0 \cap \mathcal{L}_{E}) \preceq \lexmax(\mathcal{L}_{0}\cap\mathcal{X}_{\mathcal{S}_{e^*}}\cap\mathcal{A}_{\Pi}).
\end{align*}
To complete the proof, note that for any $\Lambda\in\mathcal{L}_{0}\cap\mathcal{X}_{\mathcal{S}_{e^*}}\cap\mathcal{A}_{\Pi}$, we have $e^*\in{\bf E}(W+P)$ for all $P\in\Lambda$ by Lemma~\ref{lem:exclusive}. Therefore,
\begin{align*}
\lexmax(\mathcal{L}_{0}\cap\mathcal{X}_{\mathcal{S}_{e^*}}\cap\mathcal{A}_{\Pi}) \preceq  \lexmax(\mathcal{L}_0 \cap \mathcal{L}_{E}),
\end{align*}
which completes the proof.
\end{proof}

\begin{remark}
In Theorem~\ref{thm:eReformulation}, we note that~\eqref{eq:lSet4Pi}-\eqref{eq:lSet4Exclusive} are defined only by the bounds $\underline{\lambda}_{e}$ for all $e\in\mathcal{S}_{e^*}\cup\{e^*\}$ and $\overline{\lambda}_{e}$ for all $e\in\Pi$. All other upper and lower bounds are not involved in the set definitions, and are thus unbounded in the lexicographic maximisation. We use $\overline{\lambda}_{e} = \infty$ to denote that any positive perturbation is allowable, similarly using $-\underline{\lambda}_{e} = -\infty$ for negative perturbations.
\end{remark}

In Theorem~\ref{thm:eReformulation} we reformulate Problem~\ref{prob:ESensitivity}, using auxiliary variables $\mathcal{S}_{e^*}\in{\bf S}_{e^*}$ and $\Pi\in\mathcal{P}_{e^*}$. In the following section we construct $\Lambda^* = \lexmax(\mathcal{L}_{0}\cap\mathcal{X}_{\mathcal{S}_{e^*}}\cap\mathcal{A}_{\Pi})$ for a given exclusive set $\mathcal{S}_{e^*}$ and assignment $\Pi$. Finally, we conclude with algorithms to construct the particular exclusive set $\widehat{\mathcal{S}}_{e^*}\in{\bf S}_{e^*}$ and assignment $\widehat{\Pi}\in\mathcal{P}_{e^*}$ such that
\begin{align*}
{ \forall \mathcal{S}_{e^*}\in{\bf S}_{e^*}, \; \forall \Pi\in\mathcal{P}_{e^*} :\;}\lexmax(\mathcal{L}_{0}\cap\mathcal{X}_{\widehat{S}_{e^*}}\cap\mathcal{A}_{\widehat{\Pi}}) \succeq  \lexmax(\mathcal{L}_{0}\cap\mathcal{X}_{\mathcal{S}_{e^*}}\cap\mathcal{A}_{\Pi}),
\end{align*}
completing the solution to Problem~\ref{prob:ESensitivity}.

\subsection{Constructing Edge Allowable Perturbation Intervals}\label{subsec:eLambdaFromSnPi}

In this section we construct $\lexmax(\mathcal{L}_{0}\cap\mathcal{X}_{\mathcal{S}_{e^*}}\cap\mathcal{A}_{\Pi})$ for a given exclusive set $\mathcal{S}_{e^*}$ and assignment $\Pi$. We first define two bounds based on these sets
\begin{align}
\overline{\lambda}_{e^*} := \min\limits_{e\in\mathcal{S}_{e^*}} \frac{w_{e}-w_{e^*}}{2}\;, \qquad \underline{\lambda}_{e^*} := \min\limits_{e\in\Pi\setminus e^*} \frac{w_{e^*}-w_{e}}{2}. \label{eq:bottleneckPerturbations}
\end{align}
These bounds represent the weight gap between the bottleneck edge and the ``closest'' edges from both $\mathcal{S}_{e^*}$ and $\Pi$, where we take closest here to mean the minimum weight difference. With these bounds, we construct an array of intervals $\Lambda := ([-\underline{\lambda}_e,\overline{\lambda}_e])_{\forall e\in\mathcal{E}}$ such that
\begin{align}
\forall e\in\mathcal{E}\;:[-\underline{\lambda}_e,\overline{\lambda}_e]=
\begin{cases}
[w_{e^*}+\overline{\lambda}_{e^*}-w_{e},\infty] \; &\textrm{if} \; e\in\mathcal{S}_{e^*}\\
[-\infty,w_{e^*}-\underline{\lambda}_{e^*}-w_{e}] \; &\textrm{if} \; e\in\Pi\setminus e^*\\
[-\underline{\lambda}_{e^*},\overline{\lambda}_{e^*}]\; &\textrm{if} \; e = e^* \\
[-\infty,\infty]\; &\textrm{if} \; e\not\in\mathcal{S}_{e^*}\cup\Pi
\end{cases}.
\label{eq:givenSets}
\end{align} 

\begin{lemma}\label{lem:givenSetsProof}
For a bipartite graph $\mathcal{G}$ with weights $W$, let $\Pi\in\mathcal{P}_{e^*}$ and $\mathcal{S}_{e^*}$ be a given assignment and exclusive set satisfying~\eqref{eq:maxWeightEdge} and \eqref{eq:lessThanExclusive}, respectively, relative to $e^*\in{\bf E}(W)$. Then $\Lambda$ as defined in~\eqref{eq:givenSets} satisfies
\begin{align*}
\Lambda = \lexmax (\mathcal{L}_{0}\cap\mathcal{X}_{\mathcal{S}_{e^*}}\cap\mathcal{A}_{\Pi})
\end{align*}
where $\mathcal{A}_{\Pi}$ and $\mathcal{X}_{\mathcal{S}_{e^*}}$ are given by \eqref{eq:lSet4Pi}-\eqref{eq:lSet4Exclusive}.
\end{lemma}

\begin{proof}
The exclusive set $\mathcal{S}_{e^*}$ is used to construct a set of constraints in~\eqref{eq:lSet4Exclusive} of the form
\begin{align}
{ \forall e\in\mathcal{S}_{e^*}:\;}w_{e^*} + \overline{\lambda}_{e^*} \leq w_{e} - \underline{\lambda}_{e}. \label{eq:exclConstraint}
\end{align}
Note that the value $\overline{\lambda}_{e^*}$ in~\eqref{eq:bottleneckPerturbations} is less than or equal to the bounds $\underline{\lambda}_{e}$ for all $e\in\mathcal{S}_{e^*}$ in~\eqref{eq:givenSets}. If all of the bounds were larger than $\overline{\lambda}_{e^*}$, then there would exist an edge $e\in\mathcal{S}_{e^*}$ such that ~\eqref{eq:exclConstraint} is violated. Therefore, for $\Lambda$ to be the lexicographic max, we have $\overline{\lambda}_{e^*}$ must be as defined in~\eqref{eq:bottleneckPerturbations}. Then, given that $\overline{\lambda}_{e^*}$ is determined, the bounds $\underline{\lambda}_{e}$ are trivially maximised by solving~\eqref{eq:exclConstraint}, resulting in the bounds from~\eqref{eq:givenSets}.

For the edges in $\Pi\setminus e^*$, the proof is similar. The edges in $\Pi\setminus e^*$ are used to construct a set of constraints in~\eqref{eq:lSet4Pi}
\begin{align}
\forall e\in\Pi\setminus e^*:\;w_{e^*}-\underline{\lambda}_{e^*} \geq w_{e} + \overline{\lambda}_{e}. \label{eq:piConstraint}
\end{align}
The value of $\underline{\lambda}_{e^*}$ in~\eqref{eq:bottleneckPerturbations} is less than or equal to the bounds $\overline{\lambda}_{e}$ for all $e\in\Pi\setminus e^*$ in~\eqref{eq:givenSets}. If all of the bounds were larger than $\underline{\lambda}_{e^*}$, then there would exist an edge $e\in\Pi\setminus e^*$ such that~\eqref{eq:piConstraint} is violated. Therefore, for $\Lambda$ to be the lexicographically max, we have $\underline{\lambda}_{e^*}$ must be as defined in~\eqref{eq:bottleneckPerturbations}. Then, given that $\underline{\lambda}_{e^*}$ is determined, the bounds $\overline{\lambda}_{e}$ are trivially maximised by solving~\eqref{eq:piConstraint}, resulting in the bounds from~\eqref{eq:givenSets}. 

We have shown that the lower bounds $\underline{\lambda}_e$ for all $e\in\mathcal{S}_{e^*}\cup\{e^*\}$ and the upper bounds $\overline{\lambda}_e$ for all $e\in\Pi$ are lexicographically maximised. The remaining bounds are set to $\infty$, so for a given assignment $\Pi$ and exclusive set $\mathcal{S}_{e^*}$ we have $\Lambda = \lexmax (\mathcal{L}_{0}\cap\mathcal{X}_{\mathcal{S}_{e^*}}\cap\mathcal{A}_{\Pi})$ with $\Lambda$ as constructed in~\eqref{eq:givenSets}.
\end{proof}

\begin{remark}
For a bipartite graph $\mathcal{G}=(\mathcal{V},\mathcal{E})$, constructing $\Lambda$ as defined in~\eqref{eq:givenSets} has computational complexity $\mathcal{O}(|\mathcal{E}|)$.
\end{remark}

\subsection{Constructing the Exclusive Set and Assignment}\label{subsec:eOptimalSets}

We now show how to select the exclusive set $\widehat{\mathcal{S}}_{e^*}$ and bottleneck assignment $\widehat{\Pi}$ such that
\begin{align*}
{ \forall \mathcal{S}_{e^*}\in{\bf S}_{e^*},\; \forall \Pi\in\mathcal{P}_{e^*}:\;} \lexmax(\mathcal{L}_{0}\cap\mathcal{X}_{\widehat{\mathcal{S}}_{e^*}}\cap\mathcal{A}_{\widehat{\Pi}}) \succeq  \lexmax(\mathcal{L}_{0}\cap\mathcal{X}_{\mathcal{S}_{e^*}}\cap\mathcal{A}_{\Pi}).
\end{align*}
We begin with the bottleneck assignment $\widehat{\Pi}$.

\begin{definition}[Lexicographic Assignment]\label{def:lexAssignment}
For weight matrix $W$, an assignment $\Pi$ is a lexicographic assignment if there does \textbf{not} exist a $\Pi'$ such that
\begin{align*}
\phi_i(\Pi',W) < \phi_i(\Pi,W) \quad \textrm{and} \quad \forall j<i:\;\phi_j(\Pi',W) = \phi_j(\Pi,W),
\end{align*}
where $\phi_i(\Pi,W)$ returns the $i$-th largest element of the set $\{w_{e}| e\in\Pi\}$, i.e. the $i$-th largest weight.
\end{definition}

Informally, the assignment in which the largest weight is minimised, the second largest weight is minimised, \dots is the lexicographic assignment. The lexicographic assignment can be seen as a continuation of the bottleneck assignment objective, and can be computed in polynomial time as in \cite{burkard1991lexicographic}.



We now turn to the choice of the exclusive set $\widehat{\mathcal{S}}_{e^*}$ which corresponds to the lexicographically largest array of intervals. Note in Algorithm~\ref{alg:exclusiveSet} we use $\mathcal{S}_{e^*} \leftarrow \mathcal{S}_{e^*}\cup \{e\}$ to denote the addition of the edge $e$ to the set of edges $\mathcal{S}_{e^*}$.

\begin{algorithm}
\caption{Algorithm to Construct $\widehat{\mathcal{S}}_{e^*}$ \label{alg:exclusiveSet}}
\KwData{$W,e^*$}
\KwResult{$\widehat{\mathcal{S}}_{e^*}$}
$\widehat{\mathcal{S}}_{e^*} \leftarrow \{\}$ \;
$w_{e^*} \leftarrow \infty$ ; \phantom{12345678}\tcp{$\infty$ may be any suitably large constant}
$e \leftarrow {\bf E}(W)$ \;
\While{$w_e \neq \infty$}{
    $\widehat{\mathcal{S}}_{e} \leftarrow \widehat{\mathcal{S}}_{e^*}\cup \{e\}$ ; \tcp{Add edge $e$ to $\mathcal{S}_{e^*}$}
    $w_e \leftarrow \infty$ ; \phantom{123456}\tcp{Remove $e$ from the graph}
    $e \leftarrow {\bf E}(W)$ ; \phantom{1234}\tcp{Find the new bottleneck edge}
}
return $\widehat{\mathcal{S}}_{e^*}$
\end{algorithm}

In Algorithm~\ref{alg:exclusiveSet}, the initial bottleneck edge is found, removed from the graph by setting its weight to $\infty$, and the subsequent bottleneck edge is found. If the new bottleneck edge has finite weight, it is added to the set $\widehat{\mathcal{S}}_{e^*}$. This repeats until the bottleneck edge has weight $\infty$, i.e. there are no assignments left with the remaining edges of the graph. If, rather than setting the weight to $\infty$, the edge is removed from the graph at each iteration, then the algorithm terminates when no solution exists. We require the following assumption.

\begin{assumption}\label{ass:uniqueBottleneck}
For the bipartite graph $\mathcal{G}$ with weights $W$, the bottleneck edge $e\in{\bf E}(W)$ at each iteration of Algorithm~\ref{alg:exclusiveSet} is unique.
\end{assumption}

Assumption~\ref{ass:uniqueBottleneck} is trivially satisfied if all edge weights are distinct. If Assumption~\ref{ass:uniqueBottleneck} does not hold, Algorithm~\ref{alg:exclusiveSet} will still return an exclusive set, but we cannot certify that the resulting $\Lambda$ from~\eqref{eq:givenSets} will be the lexicographically largest array of allowable perturbations.

\begin{remark}
Assumption~\ref{ass:uniqueBottleneck}\footnote{This remark applies to  Assumption~\ref{ass:bapUniqueBottleneck} as stated later in the paper as well.} is necessary to certify the lexicographic optimality of the resulting intervals. The non-uniqueness of the choice of bottleneck edge at an iteration can affect the set of subsequent bottleneck edges. We provide Example~\ref{ex:ambiguousBottleneck} to demonstrate this effect.
\end{remark}

\begin{example}\label{ex:ambiguousBottleneck}
Consider the following weight matrix
\begin{align*} W =
\begin{bmatrix}
 0 & 10 & 0 \\
 100 & 1 & 5 \\
 0 & 5 & 0
\end{bmatrix}.
\end{align*}
The \emph{unique} bottleneck edge of this weight matrix is $e^*=(2,2)$, with weight $1$. The lexicographic assignment is $\Pi = \{(1,1),(2,2),(3,3)\}$. Algorithm~\ref{alg:exclusiveSet} begins by removing $e^*$ from the graph, yielding
\begin{align*} W =
\begin{bmatrix}
 0 & 10 & 0 \\
 100 & \infty & 5 \\
 0 & 5 & 0
\end{bmatrix},
\end{align*}
in which there are two bottleneck edges $e_1 = (2,3)$ and $e_2 = (3,2)$ with weight $5$. Let cases $1,2$ be the removal of $e_1$ and $e_2$ in subsequent iterations respectively. The exclusive sets computed for case $1,2$ are $\mathcal{S}^{1}_{e} = \{(2,3),(2,1)\}$ and $\mathcal{S}^{2}_{e} = \{(3,2),(1,2)\}$. The intervals $\Lambda_1$ and $\Lambda_2$ from~\eqref{eq:givenSets} are shown in Table~\ref{tble:cases}.
\begin{table}[ht]
\begin{center}
\caption{Ambiguous Bottleneck: $\Lambda_1$ on left, $\Lambda_2$ on right \label{tble:cases}}
\begin{tabular}{c  c  c |}
$(-\infty, 0.5]$ & $(-\infty, \infty)$ & $(-\infty, \infty)$ \\
$[-97,\infty)$ & $[-0.5, 2]$ & $[-2, \infty)$ \\
$(-\infty, \infty)$ & $(-\infty,\infty)$ & $(-\infty, 0.5]$
\end{tabular}
\begin{tabular}{c  c  c }
$(-\infty, 0.5]$ & $[-7, \infty)$ & $(-\infty, \infty)$ \\
$(-\infty, \infty)$ & $[-0.5, 2]$ & $(-\infty, \infty)$ \\
$(-\infty, \infty)$ & $[-2,\infty)$ & $(-\infty, 0.5]$
\end{tabular}
\end{center}
\end{table}
As is shown in Table~\ref{tble:cases}, the choice between removing $e_1,e_2$ introduced an identical bound of magnitude $2$. However, the choice also changed which of edges $(2,1)$ or $(1,2)$ became the subsequent bottleneck edge, which lead to the array of intervals $\Lambda_1$ being lexicographically larger than $\Lambda_2$.
\end{example}

\begin{theorem}\label{thm:eSens}
Given a weight matrix $W$ with a bottleneck edge $e^*\in{\bf E}(W)$, let $\widehat{\mathcal{S}}_{e^*},\widehat{\Pi}$ be the set of edge returned from Algorithm~\ref{alg:exclusiveSet} and the lexicographic assignment as defined in Definition~\ref{def:lexAssignment} respectively. Then $\widehat{\mathcal{S}}_{e^*}$ is an exclusive set relative to $e^*$, and if Assumption~\ref{ass:uniqueBottleneck} holds, we have
\begin{align*}
\forall \Pi\in\mathcal{P}_{e^*},\; \forall \mathcal{S}_{e^*}\in{\bf S}_{e^*}:\;\lexmax(\mathcal{L}_{0}\cap\mathcal{X}_{\widehat{\mathcal{S}}_{e^*}}\cap\mathcal{A}_{\widehat{\Pi}}) \succeq  \lexmax(\mathcal{L}_{0}\cap\mathcal{X}_{\mathcal{S}_{e^*}}\cap\mathcal{A}_{\Pi}).
\end{align*}
Therefore, the array of intervals $\Lambda = \lexmax(\mathcal{L}_{0}\cap\mathcal{X}_{\widehat{\mathcal{S}}_{e^*}}\cap\mathcal{A}_{\widehat{\Pi}})$ defined as in~\eqref{eq:givenSets} relative to $\widehat{\mathcal{S}}_{e^*},\widehat{\Pi}$ solves Problem~\ref{prob:ESensitivity}.
\end{theorem}

\begin{proof}
As noted in the proof of Lemma~\ref{lem:givenSetsProof}, the bounds derived from $\widehat{\mathcal{S}}_{e^*}$ and $\widehat{\Pi}$ are entirely independent in~\eqref{eq:givenSets}. We begin with the analysis of the set $\widehat{\mathcal{S}}_{e^*}$, and finish with $\widehat{\Pi}$. Let $\Lambda'=\lexmax(\mathcal{L}_{0}\cap\mathcal{L}_{E})$, for $\mathcal{L}_{0},\mathcal{L}_{E}$ defined in~\eqref{eq:L0} and~\eqref{eq:eAllowableSet} respectively. From Theorem~\ref{thm:eReformulation}, there exist $\mathcal{S}'_{e^*},\Pi'$ such that $\Lambda'=\lexmax(\mathcal{L}_{0}\cap\mathcal{X}_{\mathcal{S}'_{e^*}}\cap\mathcal{A}_{\Pi'})$. We will show that $\mathcal{S}'_{e^*}$ must be equal to $\widehat{\mathcal{S}}_{e^*}$ and $\Pi'$ must be equal to $\widehat{\Pi}$.

We first show $\widehat{\mathcal{S}}_{e^*}$ is an exclusive set relative to $e^*$, as defined in Definition~\ref{def:exclusiveSet}. The set $\widehat{\mathcal{S}}_{e^*}$ satisfies the first two properties from the definition trivially. To see that $\widehat{\mathcal{S}}_{e^*}$ satisfies the final property, note that each edge in the set $\widehat{\mathcal{S}}_{e^*}$, as well as $e^*$, is assigned the weight $\infty$ in the course of Algorithm~\ref{alg:exclusiveSet}. The algorithm terminates when the bottleneck edge of the graph has weight $\infty$. Therefore, there are no assignments remaining in the graph which do not assign either the bottleneck edge $e^*$ or an edge in $e\in\widehat{\mathcal{S}}_{e^*}$, which is the final property.

To show that $\mathcal{S}'_{e^*} = \widehat{\mathcal{S}}_{e^*}$, let $\{\hat{e}^{(1)},\hat{e}^{(2)},...\}$ be the edges in $\widehat{\mathcal{S}}_{e^*}$ ordered by weight from least to greatest, i.e. $w_{\hat{e}^{(k)}} \leq w_{\hat{e}^{(k+1)}}$ for all $k\in[1,2,...,|\widehat{\mathcal{S}}_{e^*}|]$. This ordering is also the order in which the edges are added to $\widehat{\mathcal{S}}_{e^*}$ in Algorithm~\ref{alg:exclusiveSet}. To see this, note that removing an edge from the graph can only increase the bottleneck edge weight. We similarly define $\{e'^{(1)},e'^{(2)},...\}$ to be the weight ordered edges in $\mathcal{S}'_{e^*}$. Assume there exists some $k\in[0,1,2,...,|\widehat{\mathcal{S}}_{e^*}|]$ such that $\hat{e}^{(j)} = e'^{(j)}$ for all $j\leq k$, i.e. the sets $\mathcal{S}'_{e^*},\widehat{\mathcal{S}}_{e^*}$ share the first $k$ minimum weight edges, but may differ from $k+1$ onwards. Note that $k$ can be zero. If $k+1 < |\widehat{\mathcal{S}}_{e^*}|$, then edge $\hat{e}^{(k+1)}$ is the bottleneck edge of the graph at iteration $k+1$ in Algorithm~\ref{alg:exclusiveSet}. Further, by Assumption~\ref{ass:uniqueBottleneck}, there exists an assignment $\Pi_{(k+1)}$ such that $\hat{e}^{(k+1)}$ is the unique maximum weight edge. As $\mathcal{S}'_{e^*}$ is an exclusive set, at least one edge $e\in\Pi_{(k+1)}$ must be in $\mathcal{S}'_{e^*}$. All previously found edges $\hat{e}^{(j)}$ for all $j<k$ have been removed from the graph, so $\hat{e}^{(j)}\not\in\Pi_{(k+1)}$ for all $j<k$. From~\eqref{eq:givenSets}, we see that the magnitude of the bound $\underline{\lambda}'_{e}$ for an edge $e\in\mathcal{S}'_{e}$ increases as a function of $w_{e}$, and because $\hat{e}^{(k+1)}$ is the unique maximum weight edge in $\Pi_{(k+1)}$, we must have $\hat{e}^{(k+1)}\in\mathcal{S}'_{e^*}$ and $e'^{(k+1)}=\hat{e}^{(k+1)}$. Any other choice of $e\in\Pi_{(k+1)}$ would result in a lexicographically smaller $\Lambda'$, which contradicts the assumptions that $\Lambda'=\lexmax(\mathcal{L}_{0}\cap\mathcal{L}_{E})$. This argument holds until $k+1 > |\widehat{\mathcal{S}}_{e^*}|$, at which point we have $\widehat{\mathcal{S}}_{e^*}\subseteq \mathcal{S}'_{e^*}$. As $\widehat{\mathcal{S}}_{e^*}$ is an exclusive set, adding any more edges to $\mathcal{S}'_{e^*}$ would yield a lexicographically smaller $\Lambda'$, and thus $\mathcal{S}'_{e^*} = \widehat{\mathcal{S}}_{e^*}$.

To show that $\Pi'=\widehat{\Pi}$, we note that the bounds in~\eqref{eq:givenSets} are strictly increasing as a function of $w_{e^*} - w_{e}$ for $e\in\Pi$. Therefore, to lexicographically maximise the bounds from $\widehat{\Pi}$, we must lexicographically minimise $w_{e}$ for all $e\in\widehat{\Pi}$. Therefore, by the definition of the lexicographic assignment, we have $\Pi' = \widehat{\Pi}$.

By assuming $\Lambda'=\lexmax(\mathcal{L}_{0}\cap\mathcal{L}_{E})$, we showed that $\Lambda'=\lexmax(\mathcal{L}_{0}\cap\mathcal{X}_{\mathcal{S}'_{e^*}}\cap\mathcal{A}_{\Pi'})$ by Theorem~\ref{thm:eReformulation} and that $\mathcal{S}'_{e^*} = \widehat{\mathcal{S}}_{e^*}$ and $\Pi'=\widehat{\Pi}$. Therefore, the array of intervals $\Lambda^*=\lexmax(\mathcal{L}_{0}\cap\mathcal{X}_{\widehat{\mathcal{S}}_{e^*}}\cap\mathcal{A}_{\widehat{\Pi}})$ constructed as in~\eqref{eq:givenSets} with $\widehat{\mathcal{S}}_{e^*},\widehat{\Pi}$ solves Problem~\ref{prob:ESensitivity}.
\end{proof}


We conclude this section by applying the results to a small numerical example.

\begin{example}\label{ex:edgeDemo}
Consider a scenario where the number of vertices in each set is $n=m=3$. We consider a case where all edges have finite weight, although in general this is not necessary. Consider the following weight matrix
\begin{align*} W =
\begin{bmatrix}
 2 & 91 & 63 \\
 26 & 89 & 93 \\
 48 & 60 & 71
\end{bmatrix}.
\end{align*}
The bottleneck edge of this weight matrix is $e^*=(1,3)$, with weight $63$. We define the set $\widehat{\mathcal{S}}_{e^*}$ as it is returned from Algorithm~\ref{alg:exclusiveSet}, $\widehat{\mathcal{S}}_{e^*}=\{(2,2),(1,2),(2,3)\}$.
With the lexicographic assignment $\widehat{\Pi} = \{(2,1),(3,2),(1,3)\}$ and the aforementioned $\widehat{\mathcal{S}}_{e^*}$, the lexicographically largest allowable perturbation intervals array is given by $$\Lambda^* = \lexmax(\mathcal{L}_{0}\cap\mathcal{X}_{\mathcal{S}_{e^*}}\cap\mathcal{A}_{\Pi}),$$
where $\mathcal{L}_{0}, \mathcal{X}_{\mathcal{S}_{e^*}}, \mathcal{A}_{\Pi}$
are defined as in Theorem~\ref{thm:eReformulation}. The elements of interval array $\Lambda^*$ are presented in Table~\ref{tble:edgeDemo}.

\begin{table}[ht]
\begin{center}
\caption{Lexicographically Largest Allowable Perturbation Intervals, $\Lambda^*$ \label{tble:edgeDemo}}
\begin{tabular}{c  c  c }
$(-\infty, \infty)$ & $[-15, \infty)$ & $[-1.5,13]$ \\
$(-\infty,35.5]$ & $[-13, \infty)$ & $[-17, \infty)$ \\
$(-\infty, \infty)$ & $(-\infty,1.5]$ & $(-\infty, \infty)$
\end{tabular}
\end{center}
\end{table}
It is simple to verify that any perturbation $P\in\Lambda^*$ satisfies $e^*\in{\bf E}(W+P)$.
\end{example}
\section{Bottleneck Assignment Sensitivity}\label{sec:sensBttleAss}
In this section we solve Problem~\ref{prob:BAPSensitivity} by providing Algorithm~\ref{alg:bapSens} to construct $\Lambda^*= \lexmax(\mathcal{L}_0\cap\mathcal{L}_{A})$. In the previous section, we were able to separate the construction of the exclusive set from the construction of the intervals. However, for Problem~\ref{prob:BAPSensitivity}, we must construct them simultaneously. We begin with a similar result to Theorem~\ref{thm:eReformulation}, reformulating the bottleneck assignment sensitivity problem.

\begin{theorem}\label{thm:bapReformulation}
For a given graph $\mathcal{G}$ with weights $W$, let $\Pi\in{\bf BAP}(W)$ be a bottleneck assignment. There exist $\mathcal{S}_{e}\in{\bf S}_{e}$ for all $e\in\Pi$ such that
\begin{align}
\lexmax \left( \mathcal{L}_0\cap_{e\in\Pi}\mathcal{X}_{\mathcal{S}_{e}} \right) = \lexmax \left(\mathcal{L}_{0}\cap\mathcal{L}_{A} \right) \label{eq:bapReformulation}
\end{align}
for $\mathcal{L}_0,\mathcal{X}_{\mathcal{S}_{e}},\mathcal{L}_{A}$ defined in~\eqref{eq:L0},\eqref{eq:lSet4Exclusive}, and \eqref{eq:bapAllowableSet} respectively.
\end{theorem}

\begin{proof}
This proof is similar to the proof of Theorem~\ref{thm:eReformulation}. Assume $\Lambda^* = \lexmax \left(\mathcal{L}_{0}\cap\mathcal{L}_{A} \right)$, where $\mathcal{L}_{A}$ is defined relative to a given assignment $\Pi$. We define a set of perturbations $P_{e_i}$ for all $e_i\in\Pi$ by
\begin{align}
P_{e_i}[e] :=
\begin{cases}
-\underline{\lambda}^{*}_{e} &\textrm{if} \; e \neq e_i\\
\overline{\lambda}^{*}_{e} &\textrm{if} \; e = e_i
\end{cases}
\end{align}
such that every edge perturbation is the lower bound of its corresponding interval in $\Lambda^*$, except $e_i$. By~\eqref{eq:bapAllowableSet}, we have that $\Pi\in{\bf BAP}(W+P)$ for all $P\in\Lambda^*$, and in particular, for any $P_{e_i}$. By Lemma~\ref{lem:exclusive}, there exists an exclusive set $\mathcal{S}_{e_i}$ relative to $e_i$ which satisfies~\eqref{eq:lessThanExclusive} for $W_{e_i} := W + P_{e_i}$. By the construction of $P_{e_i}$ we have $\Lambda^*\in\mathcal{X}_{\mathcal{S}_{e_i}}$ as
\begin{align*}
{\forall e\in\mathcal{S}_{e_i},\;\forall P\in\Lambda^*:\;} w_{e_i}+P[e_i] &\leq w_{e_i} + \overline{\lambda}^{*}_{e_i} \leq w_{e} - \underline{\lambda}^{*}_{e} \leq w_{e} - P[e]\\
\intertext{for all $e_i\in\Pi$. Therefore, }
\lexmax \left( \mathcal{L}_0\cap_{e\in\Pi}\mathcal{X}_{\mathcal{S}_{e}} \right) &\succeq  \lexmax \left(\mathcal{L}_{0}\cap\mathcal{L}_{A} \right) = \Lambda^*.
\end{align*}
To show the converse, consider any perturbation $P\in\mathcal{L}_0\cap_{e\in\Pi}\mathcal{X}_{\mathcal{S}_{e}}$. Let $e_i\in\Pi$ be any of the edges in the assignment with maximum weight in the perturbed weights, i.e. $e_i\in\argmax_{e_i\in\Pi} w_{e_i} + P[e_i]$. By the definition of $\mathcal{X}_{\mathcal{S}_{e_i}}$, the exclusive set $\mathcal{S}_{e_i}$ satisfies~\eqref{eq:lessThanExclusive}, and because it is a maximum weight edge of the edges in $\Pi$, we have the $\Pi$ satisfies~\eqref{eq:maxWeightEdge}. Therefore, by Lemma~\ref{lem:exclusive}, we have that $e_i$ is a bottleneck edge, and $\Pi$ is a bottleneck assignment. Therefore,
\begin{align}
\lexmax \left( \mathcal{L}_0\cap_{e\in\Pi}\mathcal{X}_{\mathcal{S}_{e}} \right) \preceq  \lexmax\left(\mathcal{L}_{0}\cap\mathcal{L}_{A} \right) = \Lambda^* ,
\end{align}
completing the proof.
\end{proof}

In Theorem~\ref{thm:bapReformulation}, we reformulate Problem~\ref{prob:BAPSensitivity} using exclusive sets, similarly to Theorem~\ref{thm:eReformulation}. Algorithm~\ref{alg:bapSens} constructs the exclusive sets $\mathcal{S}_{e}$ for all $e\in\Pi$, from Theorem~\ref{thm:bapReformulation}. Similarly to Algorithm~\ref{alg:exclusiveSet}, in Algorithm~\ref{alg:bapSens} we find the bottleneck edge of a graph and ``remove'' it, i.e. setting the weight to $\infty$. However, we must do this in parallel for each edge $e\in\Pi$. The primary difference between Algorithm~\ref{alg:bapSens} and Algorithm~\ref{alg:exclusiveSet} is the construction of a new matrix of weights $B_{e}$ with each element $B_{e}[e'] := \textrm{boundValue}(w_{e},w_{e'},\overline{\lambda}_{e},\underline{\lambda}_{e'}$ for edges $e\in\Pi$ at each iteration.

\begin{function}
\caption{boundValue()$(w_{e},w_{e'},\overline{\lambda}_{e},\underline{\lambda}_{e'})$\label{func:maxBound}}
\tcp{$\lambda = \infty$ indicates an undetermined bound}
\uIf{$\overline{\lambda}_{e} = \infty \wedge \underline{\lambda}_{e'} = \infty$}{
  return $\frac{w_{e'} - w_{e}}{2}$\;
}
\uElseIf{$\overline{\lambda}_{e} = \infty$}{
  return $ w_{e'} - \underline{\lambda}_{e'} - w_{e}$\;
}
\uElseIf{$\underline{\lambda}_{e'} = \infty$}{
  return $ w_{e'} - w_{e} - \overline{\lambda}_{e}$\;
}
\uElseIf{$w_{e} + \overline{\lambda}_{e} \leq w_{e'} - \underline{\lambda}_{e'}$}{
  return $\infty$\;
}
\Else{
  return $-\infty$\;
}
\end{function}

To motivate the definition of function boundValue, assume that at some iteration of the algorithm we have a set of determined upper and lower bounds $\overline{\lambda}_{e},\underline{\lambda}_{e}$, along with the sets of edges $\mathcal{S}_{e}$, which are not yet exclusive sets as the algorithm has not terminated. Adding a new edge $e'$ to the set $\mathcal{S}_{e}$, by the definition of~\eqref{eq:lSet4Exclusive}, introduces a new constraint of the form
\begin{align}
w_{e} + \overline{\lambda}_{e} \leq w_{e'} - \underline{\lambda}_{e'}.\label{eq:introducedConstraint}
\end{align}
Assuming the previously determined bounds are fixed, the values $\overline{\lambda}_{e},\underline{\lambda}_{e'}$ return from boundValue$(w_{e},w_{e'},\overline{\lambda}_{e},\underline{\lambda}_{e'})$ are the lexicographically largest values which satisfy the constraint~\eqref{eq:introducedConstraint}. Previously determined bounds are considered fixed as the algorithm determines them in non-decreasing order, proof in Lemma~\ref{lem:nonDecreasing}. Lowering a previously determined bound to increase a subsequent larger bound is lexicographically sub-optimal. Algorithm~\ref{alg:bapSens}, which solves Problem~\ref{prob:BAPSensitivity}, is presented below. Note that $B_e[e_j]$ is used to indicate the element of $B_e$ corresponding to edge $e_j$ in Algorithm~\ref{alg:bapSens}.

\begin{algorithm}
\caption{Constructing $\Lambda = \lexmax \left( \mathcal{L}_0\cap\mathcal{L}_{A} \right)$\label{alg:bapSens}}
\KwData{$\mathcal{G},W,\Pi$}
$\Lambda \leftarrow \{[-\underline{\lambda}_{e},\overline{\lambda}_{e}] \mid { \forall e\in\mathcal{E}:\;}\underline{\lambda}_{e}\leftarrow \infty\; ,\; \overline{\lambda}_{e}\leftarrow\infty\}$ \;
\Repeat{
  \For{$e\in\Pi$}{
    \For{$e'\in\mathcal{E}\setminus\{e\}$}{
        $B_{e}[e'] = \textrm{boundValue}(w_{e},w_{e'},\overline{\lambda}_{e},\underline{\lambda}_{e'})$ \;
    }
    $b_{e} \leftarrow {\bf E}(B_{e})$ \;
  }
  $\widehat{e} \leftarrow \argmin\limits_{e\in\Pi} B_{e}[b_{e}]$ ; \tcp{Edge $\hat{e}\in\Pi$ has tightest bound $B_{\hat{e}}[b_{\widehat{e}}]$}
  \uIf{$B_{\widehat{e}}[b_{\widehat{e}}] = \infty$}{
    return $\Lambda$ ; \phantom{12345}\tcp{Algorithm terminates}
  }
  \uIf{$\overline{\lambda}_{\widehat{e}} = \infty$}{
    $\overline{\lambda}_{\widehat{e}} \leftarrow B_{\widehat{e}}[b_{\widehat{e}}]$ ; \phantom{12'}\tcp{Bound $\overline{\lambda}_{\widehat{e}}$ has not been previously set}
  }
  \uIf{$\underline{\lambda}_{b_{\widehat{e}}} = \infty$}{
    $\underline{\lambda}_{b_{\widehat{e}}} \leftarrow B_{\widehat{e}}[b_{\widehat{e}}]$ ; \phantom{12}\tcp{Bound $\underline{\lambda}_{b_{\widehat{e}}}$ has not been previously set}
  }
}
\end{algorithm}

We make an assumption on the structure of the weighted graph, analogous to Assumption~\ref{ass:uniqueBottleneck}.

\begin{assumption}\label{ass:bapUniqueBottleneck}
The bottleneck edges for all constructed weights in Algorithm~\ref{alg:bapSens} are unique.
\end{assumption}

Algorithm~\ref{alg:bapSens} suffers from the same lack of tie-breaking mechanism as Algorithm~\ref{alg:exclusiveSet}. For an example where Assumption~\ref{ass:uniqueBottleneck} does not hold and Algorithm~\ref{alg:exclusiveSet} is not guaranteed to work, see Example~\ref{ex:ambiguousBottleneck}. Before proving that Algorithm~\ref{alg:bapSens} solves Problem~\ref{prob:BAPSensitivity}, we establish the following lemma.

\begin{lemma}\label{lem:nonDecreasing}
Let $B_{e}^{(k)}$ be the weights computed with respect to edge $e\in\Pi$ on iteration $k$ in line $4$ of Algorithm~\ref{alg:bapSens}, and let $b^{(k)}_e = {\bf E}(B_{e}^{(k)})$ be the bottleneck edge of these weights. Then the sequence $\{\min\limits_{e\in\Pi}B_{e}^{(k)}[b^{(k)}_e]\}$ is non-decreasing.
\end{lemma}

\begin{proof}
See \ref{app:nonDecreasing}.
\end{proof}

Lemma~\ref{lem:nonDecreasing} enables the proof of Theorem~\ref{thm:bapSens} to work similarly to the proof of Theorem~\ref{thm:eSens}.

\begin{theorem}\label{thm:bapSens}
For a bipartite graph $\mathcal{G}$ with weights $W$, let $\Pi\in{\bf BAP}(W)$. If Assumption~\ref{ass:bapUniqueBottleneck} holds, then the array of intervals $\Lambda^*$ returned from Algorithm~\ref{alg:bapSens} satisfies
\begin{align*}
\Lambda^* = \lexmax \left( \mathcal{L}_0\cap\mathcal{L}_{A} \right),
\end{align*}
for $\mathcal{L}_{A}$ defined in~\eqref{eq:bapAllowableSet}. Equivalently, $\Lambda^*$ is the solution to Problem~\ref{prob:BAPSensitivity}.
\end{theorem}

\begin{proof}
We show that the array of intervals is lexicographically largest by proving that at each iteration if the magnitude of the bound was larger then the array of intervals would be unallowable.

Let $\Lambda'\in \lexmax(\mathcal{L}_0\cap\mathcal{L}_{A})$ for the given assignment $\Pi$. By Theorem~\ref{thm:bapReformulation} there exists a set of exclusive sets $\{\mathcal{S}'_{e}\}_{e\in\Pi}$ such that $\Lambda'=\lexmax \left( \mathcal{L}_0\cap_{e\in\Pi}\mathcal{X}_{\mathcal{S}'_{e}} \right)$. We will show that $\Lambda' = \Lambda^*$ for $\Lambda^*$ the intervals returned from Algorithm~\ref{alg:bapSens}.

For iteration $k$ of Algorithm~\ref{alg:bapSens} let $\hat{e}^{(k)} = \argmin\limits_{e\in\Pi}B_{e}^{(k)}[b^{(k)}_e]$ represent the edge $e\in\Pi$ for which the corresponding bottleneck edge $b^{(k)}_{\hat{e}}$ has minimum weight in $B_{e}^{(k)}$. These are the edges determined in Line $6$ of Algorithm~\ref{alg:bapSens}. During iteration $k$, if the algorithm does not terminate, then one or both of $\overline{\lambda}^*_{\hat{e}^{(k)}},\underline{\lambda}^*_{b^{(k)}_{\hat{e}}}$ will be set. After iteration $k$, both of $\overline{\lambda}^*_{\hat{e}^{(k)}},\underline{\lambda}^*_{b^{(k)}_{\hat{e}}}$ will be set and will satisfy
\begin{align}
w_{\hat{e}^{(k)}} + \overline{\lambda}^*_{\hat{e}^{(k)}}= w_{b^{(k)}_{\hat{e}}} - \underline{\lambda}^*_{b^{(k)}_{\hat{e}}}.
\end{align}
We therefore have that in subsequent iterations, the edge $b^{(k)}_{\hat{e}}$ is removed from $B_{\hat{e}^{(k)}}$ or equivalently has weight $\infty$ (Line 8 in function boundValue).

Assume there exists some $k\geq0$ such that for all $j\leq k$ we have $\overline{\lambda}^*_{\hat{e}^{(j)}} = \overline{\lambda}'_{\hat{e}^{(j)}}$ and $\underline{\lambda}^*_{b^{(j)}_{\hat{e}}} = \underline{\lambda}'_{b^{(j)}_{\hat{e}}}$, i.e. $\Lambda^*$ from Algorithm~\ref{alg:bapSens} and $\Lambda'\in \lexmax(\mathcal{L}_0\cap\mathcal{L}_{A})$ have the same $k$ minimal bound values. Note that $k$ can be zero. Then, if edge $b^{(k+1)}_{\hat{e}}$ has finite weight, it is the bottleneck edge of the graph with weights $B_{\hat{e}^{(k+1)}}$ at iteration $k+1$. Further, by Assumption~\ref{ass:bapUniqueBottleneck}, there exists an assignment $\Pi^{(k+1)}$ such that $b^{(k+1)}_{\hat{e}}$ is the unique maximum weight edge over $B_{\hat{e}^{(k+1)}}$. By Theorem~\ref{thm:bapReformulation}, $\mathcal{S}'_{\hat{e}^{(k+1)}}$ is an exclusive set, and thus at least one of the edges in $\Pi^{(k+1)}$ must be in $\mathcal{S}'_{\hat{e}^{(k+1)}}$ and satisfy
\begin{align}
w_{\hat{e}^{(k+1)}} + \overline{\lambda}'_{\hat{e}^{(k+1)}} \leq w_{b^{(k+1)}_{\hat{e}}} - \underline{\lambda}'_{b^{(k+1)}_{\hat{e}}}.
\end{align}
By Lemma~\ref{lem:nonDecreasing}, all subsequent bounds in $\Lambda^*$ will be equal to or larger than the weight of the bottleneck edge $b^{(k+1)}_{\hat{e}}$. Therefore, as $\Lambda'$ is the lexicographically largest, we must have the same property. However, as $b^{(k+1)}_{\hat{e}}$ is the unique maximum over $B_{\hat{e}^{(k+1)}}$, we must have that $\Lambda'$ shares the same bound as $\Lambda^*$ from the $k+1$ iteration, i.e. $\overline{\lambda}^*_{\hat{e}^{(k+1)}} = \overline{\lambda}'_{\hat{e}^{(k+1)}}$ and $\underline{\lambda}^*_{b^{(k+1)}_{\hat{e}}} = \underline{\lambda}'_{b^{(k+1)}_{\hat{e}}}$.

Eventually, Algorithm~\ref{alg:bapSens} reaches an iteration $l$ such that $b^{(l)}_{\hat{e}}$ has weight $\infty$, and the algorithm terminates. To see this, note that by Lemma~\ref{lem:nonDecreasing} the bottleneck edges have non-decreasing weights at each iteration, and at least one edge weight from one of $B_{e}$ for $e\in\Pi$ is set to $\infty$ at each iteration, so the algorithm must terminate as there are only finitely many edges. When Algorithm~\ref{alg:bapSens} terminates, for each of $e_i\in\Pi$, we may construct an exclusive set $\mathcal{S}^*_{e_i}$ such that $\Lambda^*\in\mathcal{X}_{\mathcal{S}^*_{e_i}}$ by checking for which edges boundValue$(w_{e},w_{e'},\overline{\lambda}_{e},\underline{\lambda}_{e'}) = \infty$. These sets are exclusive sets as each $B_{e_i}$ has bottleneck edge weight $\infty$, implying there exists no assignment which does not assign one of these edges. Further, we have $\Lambda^*\in\mathcal{X}_{\mathcal{S}^*_{e_i}}$ as function boundValue assigns a weight of $\infty$ (Line 8) only when $w_{e_i} + \overline{\lambda}^*_{e_i} \leq w_{e} - \underline{\lambda}^*_{e}$. Given that all bounds which are undetermined after Algorithm~\ref{alg:bapSens} terminates remain at $\infty$, and $\Lambda'$ is the lexicographically largest allowable array of intervals, we must have that $\Lambda' = \Lambda^*$.
\end{proof}

We conclude by demonstrating Algorithm~\ref{alg:bapSens} using the same example as in Section~\ref{sec:Esensitivity}, to highlight the differences between the measures of sensitivity.

\begin{example}\label{ex:assDemo}
Using the same weight matrix as in Example~\ref{ex:edgeDemo}, we perform a sensitivity analysis for assignment $\Pi = \{(2,1),(3,2),(1,3)\}$. The array of intervals $\Lambda$ returned from Algorithm~\ref{alg:bapSens} is shown in Table~\ref{tble:edgeDemo}.

\begin{table}[ht]
\begin{center}
\caption{Lexicographically Largest Allowable Perturbation Intervals \label{tble:assDemo}}
\begin{tabular}{c  c  c }
$(-\infty, \infty)$ & $[-15, \infty)$ & $(-\infty,13]$ \\
$(-\infty,50]$ & $[-13, \infty)$ & $[-17, \infty)$ \\
$(-\infty, \infty)$ & $(-\infty,16]$ & $(-\infty, \infty)$
\end{tabular}
\end{center}
\end{table}

We can see in this case that either the allowable perturbation intervals are identical to the ones in Example~\ref{ex:edgeDemo}, or strictly larger. However, this will not be true in general. Although the intervals $\Lambda$ corresponding to the bottleneck assignment sensitivity will be lexicographically larger, they may include more bounded edges than the intervals for the bottleneck edge sensitivity. For comparison, we will show the results from Theorem $3.1$ in~\cite{sotskov1995some}, over the same example.
\begin{table}[ht]
\begin{center}
\caption{Sensitivity Results using~\cite{sotskov1995some} \label{tble:sotsDemo}}
\begin{tabular}{c  c  c }
$[-13, 13]$ & $[-13, 13]$ & $[-13,13]$ \\
$[-13,13]$ & $[-13, 13]$ & $[-13, 13]$ \\
$[-13, 13]$ & $[-13,13]$ & $[-13, 13]$
\end{tabular}
\end{center}
\end{table}
\end{example}

\section{Computational Complexity}\label{sec:complexity}

In this section we will provide a brief accounting of the complexity of the algorithms provided, as well as a comparison to the results from~\cite{sotskov1995some}.

We first note that, in the computation of the bottleneck edge sensitivity as well as the bottleneck assignment sensitivity, we are iteratively running a bottleneck assignment solver on a graph, removing or modifying one edge weight at a time. Unsurprisingly, we can take advantage of the information from the previous iteration, to dramatically reduce the computational complexity. Using an augmenting path solver, such as from~\cite{punnen1994improved}, we may ``warm start'' the solver with the previously optimal bottleneck assignment, removing only the modified edge. This is similar to the concept used in~\cite{volgenant2006addendum}. In this case, the algorithm needs to run only a single augmenting path search with complexity $\mathcal{O}(|\mathcal{V}|)$, instead of finding bottleneck assignment solution with complexity $\mathcal{O}(|\mathcal{V}|\sqrt{|\mathcal{E}||\mathcal{V}|})$.

We begin with the computational complexity of the bottleneck edge sensitivity. To compute the bottleneck edge sensitivity of a given $e^*\in {\bf BAP}(W)$, we compute
\begin{itemize}
    \item The lexicographic assignment $\widehat{\Pi}$.
    \item The exclusive set as returned from Algorithm~\ref{alg:exclusiveSet}.
    \item The lexicographically largest $\Lambda^*$ defined by the previous two steps, as in~\eqref{eq:givenSets}.
\end{itemize}
For a bipartite graph $\mathcal{G} = (\mathcal{V}_1\cup\mathcal{V}_2,\mathcal{E})$, let $n = |\mathcal{V}_1|+|\mathcal{V}_2|$ and $m = |\mathcal{E}|$. Algorithm~\ref{alg:exclusiveSet} requires $\mathcal{O}(m)$ iterations of a bottleneck assignment solver, although as noted previously, each iteration is only a single augmenting path search with complexity $\mathcal{O}(n)$. The complexity of computing the intervals $\Lambda^*$ is then $\mathcal{O}(L(n,m)+nm+m)$, where $L(n,m)$ is the computational complexity of constructing the lexicographic assignment. For a simple example, we assume the graph is square and dense, i.e. $|\mathcal{V}_1| = |\mathcal{V}_2|$ and $\mathcal{E} = \mathcal{V}_1 \times \mathcal{V}_2$, and that the bottleneck edge of each subgraph of $\mathcal{G}$ is unique. Then the lexicographic assignment may be computed as $\mathcal{O}(n)$ iterations of the augmenting path solver, which yields a complexity of $\mathcal{O}(n^{2}+n^{3}+n^{2}) = \mathcal{O}(n^{3})$. This simplification shows that the computation of the exclusive set in Algorithm~\ref{alg:exclusiveSet}, which requires $\mathcal{O}(n^2)$ iterations of the bottleneck assignment solver, dominates the complexity.

For the complexity of the bottleneck assignment sensitivity, we note that there are $\mathcal{O}(m)$ iterations in Algorithm~\ref{alg:bapSens}, and at each iteration we update $\mathcal{O}(n)$ of the constructed weight matrices $B_{e}$ for all $e\in\Pi$. The complexity is then $\mathcal{O}(n B(n,m) + n^2m)$, for $B(n,m)$ the complexity of the bottleneck assignment solver over the given graph. Using the complexity of the bottleneck assignment algorithm from~\cite{punnen1994improved}, and assuming a square and dense graph as before, this yields a complexity of $\mathcal{O}(n^{3.5}+n^{4}) = \mathcal{O}(n^4)$.

In the introduction, we mentioned the most similar work found in the literature to the presented was~\cite{sotskov1995some}. In their work, they construct a ``sensitivity radius'' $\sigma(\mathcal{G})$, which we may interpret in the notation of this paper as the maximum scalar $\sigma\in\mathbb{R}^+$ such that all perturbations $P_\sigma := \{ P_{\sigma}[e]\in [-\sigma,\sigma] \mid e\in\mathcal{E}\}$ are allowable.
Letting $\Lambda_{\sigma}$ be the array of intervals with lower and upper bound $-\sigma,\sigma$, we then have that $\rho_1(\Lambda^*) = \rho_1(\Lambda_\sigma)$ and $\rho_i(\Lambda^*)\geq\rho_i(\Lambda_{\sigma})$ for all $i\in\{2,...,m\}$ for $\Lambda^*$ the array of intervals computed in Algorithm~\ref{alg:bapSens}.
The computation of the ``sensitivity radius'' requires solving a single bottleneck assignment, and so has complexity $\mathcal{O}(B(n,m))$ or $\mathcal{O}(n^{2.5})$ with the same simplifications as previously discussed.

\section{Conclusion}\label{sec:conclusion}
In this paper we propose two frameworks for assessing the sensitivity of a bottleneck assignment problem, as well as algorithms for the computation of the sensitivity within each framework. The analysis provided is driven primarily by the characterisation of ``exclusive sets'', and the connection of these sets to the bottleneck assignment problem. The combination of the theory of exclusive sets, along with recursive applications of the bottleneck assignment problem, allow for these algorithms to be run with off the shelf assignment solvers and a minimum of additional programming. The sensitivity analysis provided can be used to certify the optimality of a solution, if for example the true assignment weights can be shown to be contained by the provided intervals, or as a measure of the robustness of an assignment solution. Further research should concentrate on the expansion of this type of sensitivity analysis to other formulations of the assignment problem, as well as other formulations of the ``largest'' array of intervals besides the lexicographic ordering. Further, understanding the necessary modifications to relax Assumptions~\ref{ass:uniqueBottleneck} and \ref{ass:bapUniqueBottleneck} remains an open question. A weight scaling approach similar to~\cite{burkard1991lexicographic} seems to be a viable way forward to this end.



%
%
%

\section*{Acknowledgements}
This research was supported by the University of Melbourne Graduate Research Scholarship.

\appendix

\section{Proof of Lemma~\ref{lem:nonDecreasing}}\label{app:nonDecreasing}

We refer to $\min\limits_{e\in\Pi}B_{e}^{(k)}[b^{(k)}_e]$ as the bound determined on the $k$-th iteration of Algorithm~\ref{alg:bapSens}. To prove the sequence of these bounds is non-decreasing, we show that the sequence of bottleneck edge weights from each of these matrices $B_{e}$ is non-decreasing. Focusing on a single matrix $B_{e^*}$ for some $e^*\in\Pi$, we can separate the iterations into two cases: either the minimum weight bottleneck is from one of the other $B_{e'}$ corresponding to an edge $e'\in\Pi$, with $e'\neq e^*$, or it is from $B_{e^*}$.

We begin by assuming the minimum weight bottleneck edge is from one of the other graphs $b' = {\bf E}(B_e')$ with bottleneck edge weight $\lambda' = B_{e'}[b']$. If setting $\underline{\lambda}_{b'} = \lambda'$ increases or leaves unchanged the value of $B_{e^*}[b']$ then the bottleneck edge weight is trivially non-decreasing, so we assume that setting $\underline{\lambda}_{b'} = \lambda'$ results in a strict decrease of $B_{e^*}[b']$ in the subsequent iteration. We then have two cases, either $\overline{\lambda}_{e^*}$ has been determined prior to this iteration, or it has not. Assuming that $\overline{\lambda}_{e^*}$ has not been determined prior to this iteration, we have that the edge weight $B_{e^*}[b']$ had value $0.5(w_{e'}-w_{e^*})$ prior to updating $\underline{\lambda}_{b'}$, and has value $w_{e'} - w_{e^*} - \underline{\lambda}_{b'}$ after updating $\underline{\lambda}_{b'}$. If this is a strict decrease, we obtain $\underline{\lambda}_{b'} > 0.5(w_{e'}-w_{e^*})$, which may be restated as $B_{e'}[b'] > B_{e^*}[b']$ \emph{prior} to updating $\underline{\lambda}_{b'}$. Finally, given that $B_{e^*}[b^*] \leq B_{e'}[b']$ for $b^*={\bf E}(B_{e^*})$, we have $B_{e^*}[b^*] > B_{e^*}[b']$, i.e. the weight of edge $b'$ was less than the bottleneck edge $b^*$, and strictly decreased, therefore not changing the bottleneck edge weight in $B_{e^*}$. The second case, assuming $\overline{\lambda}_{e^*}$ was determined prior to the iteration, follows nearly identical logic. If the weight of the edge $B_{e^*}[b']$ decreases, it must have began with weight less than the bottleneck edge weight, and thus does not change the bottleneck edge weight in $B_{e^*}$.

We now assume the minimum weight bottleneck edge is in the graph $b^* = {\bf E}(B_{e^*})$ with weight $\lambda^* = B_{e^*}[b^*]$. Assume, previous to the $i$-th iteration, we have $\overline{\lambda}_{e^*} = \infty$, i.e. the upper bound on edge $e^*$ was undetermined. Then for any edge $e'$ satisfying $B_{e^*}[e'] \geq \lambda^*$, we have that $B_{e^*}[e']$ is increased or remains constant after updating. The subsequent bottleneck edge of $B_{e^*}$ is drawn from this set of edge weights, so the subsequent bottleneck edge weight is greater than or equal to $\lambda^*$. If instead we have $\overline{\lambda}_{e^*} \neq \infty$, then $B_{e^*}[b^*]$ is set to $\infty$ and nothing else changes, so the subsequent bottleneck edge weight from $B_e$ is increases or remains constant. \hfill$\square$



\bibliographystyle{elsarticle-harv}
\bibliography{lastBib.bib}





\end{document}